\documentclass{article}

\usepackage{times}
\usepackage{graphicx} 
\usepackage{subfigure} 
\usepackage{amsmath}
\usepackage{amssymb}
\usepackage{natbib}
\usepackage{amsthm}
\usepackage{comment}

\usepackage{algorithm}
\usepackage{algorithmic}

\newtheorem{thm}{Theorem}[section]
\newtheorem{lemma}[thm]{Lemma}
\newtheorem{corollary}[thm]{Corollary}
\newtheorem{proposition}[thm]{Proposition}

\DeclareMathOperator*{\dx}{\mathrm{d}\mathit{x}} 
\DeclareMathOperator*{\dz}{\mathrm{d}\mathit{z}} 
 
\DeclareMathOperator*{\dk}{\mathrm{d}\mathit{k}} 
\newcommand*{\LargerCdot}{\raisebox{-0.5ex}{\scalebox{1.5}{$\cdot$}}}
\newcommand*{\SmallerCdot}{\raisebox{-0.5ex}{\scalebox{1.2}{$\cdot$}}}

\usepackage{hyperref}

\usepackage[accepted]{icml2016}

\icmltitlerunning{On collapsed representation of hierarchical Completely Random Measures}

\begin{document} 

\twocolumn[
\icmltitle{On collapsed representation of hierarchical Completely Random Measures}

\icmlauthor{Gaurav Pandey}{gp88@csa.iisc.ernet.in}
\icmlauthor{Ambedkar Dukkipati}{ad@csa.iisc.ernet.in}
\icmladdress{Department of Computer Science and Automation \\
Indian Institute of Science, Bangalore-560012, India}

\icmlkeywords{boring formatting information, machine learning, ICML}

\vskip 0.3in
]
\begin{abstract}
The aim of the paper is to provide an exact approach for generating a Poisson process sampled from a hierarchical CRM, without having to instantiate the infinitely many atoms of the random measures.
We use completely random measures~(CRM) and hierarchical CRM to define a prior for Poisson processes. We derive the marginal distribution of the resultant point process, when the underlying CRM is marginalized out. Using well known properties unique to Poisson processes, we were able to derive an exact approach for instantiating a Poisson process with a hierarchical CRM prior. Furthermore, we derive Gibbs sampling strategies for hierarchical CRM models based on  Chinese restaurant franchise sampling scheme. As an example, we present the sum of generalized gamma process (SGGP), and show its application in topic-modelling. We show that one can determine the power-law behaviour of the topics and words in a Bayesian fashion, by defining a prior on the parameters of SGGP.
\end{abstract} 

\section{Introduction}

Mixed membership modelling is the problem of assigning an object to multiple latent classes/features simultaneously. Depending upon the problem, one can allow a single latent feature to be exhibited single or multiple times by the object. For instance, a document may comprise several topics, with each topic occurring in the document with variable multiplicity. The corresponding problem of mapping the words of a document to topics, is referred to as topic modelling.

While parametric solutions to mixed membership modelling have been available in literature since more than a decade~\cite{landauer1997solution, hofmann1999probabilistic, blei2001latent}, the first non-parametric approach, that allowed the number of latent classes to be determined as well, was the hierarchical Dirichlet process~(HDP)~\cite{teh2006hierarchical}. Both the approaches model the object as a set of repeated draws from an object-specific distribution, whereby the object specific distribution is itself sampled from a common distribution. On the other hand, recent approaches such as hierarchical beta-negative binomial process~\cite{zhou2011beta, broderick2011combinatorial} and hierarchical gamma-Poisson process~\cite{titsias2008infinite, zhou2013GammaPoisson} model the object as a point process, sampled from an object specific random measure, which is itself sampled from a common random measure. In some sense, these approaches are more natural for mixed membership modelling, since they model the object as a single entity rather than as a sequence of draws from a distribution.

A straightforward implementation of any of the above non-parametric models would require sampling the atoms in the non-parametric distribution for the base as well as object-specific measure. However, since the number of atoms in these distributions are often infinite, a truncation step is required to ensure tractability. 
Alternatively, for the HDP, a Chinese restaurant franchise scheme~\cite{teh2006hierarchical} can be used for collapsed inference in the model (that is, without explicitly instantiating the atoms). Fully collapsed inference scheme has also been proposed for beta-negative binomial process (BNBP)~\cite{heaukulani2013combinatorial, zhou2014beta} and Gamma-Gamma-Poisson process~\cite{zhou2015priors}. Of particular relevance is the work by~\citet{roy2014continuum}, whereby a Chinese restaurant fanchise scheme has been proposed for hierarchies of beta proceses (and its generalizations), when coupled with Bernoulli process.

In this paper, it is our aim to extend fully collapsed sampling so as to allow any completely random measure (CRM) for the choice of base and object-specific measure. As proposed in~\citet{roy2014continuum} for hierarchies of generalized beta processes, we propose Chinese restaurant franchise schemes for hierarchies of CRMs, when coupled with Poisson process. We hope that this will encourage the use of hierarchical random measures, other than HDP and BNBP, for mixed-membership modelling and will lead to further research into an understanding of the applicability of the various random measures. To give an idea about the flexibility that can be obtained by using other measures, we propose the sum of generalized gamma process (SGGP), which allows one to determine the power term in the power-law distribution of topics with documents, by defining a prior on the parameters of SGGP. Alternatively, one can also define a prior directly on the discount parameter.

The main contributions in this paper are as follows:
\begin{itemize}
\item We derive marginal distributions of Poisson process, when coupled with CRMs,
\item We provide an exact approach for generating a Poisson process sampled from a hierarchical CRM, without having to instantiate the infinitely many atoms of the random measure.
\item We provide a Gibbs sampling approach for sampling a Poisson process from a hierarchical CRM.
\item In the experiments section, we propose the sum of generalized gamma process (SGGP), and show its applicability for topic-modelling. By defining a prior on the parameters of SGGP, one can determine the power-law distribution of the topics and words in a Bayesian fashion.
\end{itemize}

\section{Preliminaries and background}
In this section, we fix the notation and recall a few well known results from the theory of point processes.
\subsection{Poisson process}
Let $(S, \mathcal{S})$ be a measurable space and $\Pi$ be a random countable collection of points on $S$. Let $N(A) = |\Pi \cap A|$, for any measurable set $A$. $N$ is also known as the counting process of $\Pi$. $\Pi$ is called a Poisson process if $N(A)$ is independent of $N(B)$, whenever $A$ and $B$ are disjoint measurable sets, and $N(A)$ is Poisson distributed with mean $\mu(A)$ for a fixed $\sigma$-finite measure $\mu$. In sequel, we refer to both the random collection $\Pi$ and its counting process $N$ as Poisson process. 

Let $(T, \mathcal{T})$ be another measurable space and $f:S\to T$ be a measurable function. If the push forward measure of $\mu$ via $f$, that is, $\mu\circ f^{-1}$ is non-atomic, then $f(\Pi) =\{f(x): x\in \Pi\}$ is also a Poisson process with mean measure $\mu\circ f^{-1}$. This is also known as the mapping proposition for Poisson processes~\cite{kingman1992poisson}. Moreover, if $\Pi_1, \Pi_2,\dotsc$ is a countable collection of independent Poisson processes with mean measures $\mu_1, \mu_2,\dotsc$ respectively, then the union $\Pi = \cup_{i=1}^\infty \Pi_i$ is also a Poisson process with mean measure $\mu = \sum_{i=1}^\infty \mu_i$. This is known as the superposition proposition. Equivalently, if $N_i$ is the counting process of $\Pi_i$, then $N = \sum_{i=1}^\infty N_i$ is the counting process of a Poisson process with mean measure $\mu = \sum_{i=1}^\infty \mu_i$.

Finally, let $g$ be a measurable function from $S$ to $\mathbb{R}$, and $\Sigma = \sum_{x\in \Pi}g(x)$. By Campbell's proposition~\cite{kingman1992poisson}, $\Sigma$ is absolutely convergent with probability, if and only if
\begin{equation}
\int_S \min(|g(x)|,1)\mu(\dx) < \infty.
\end{equation}
If this condition holds, then for any $t>0$, 
\begin{equation}
\mathbb{E}[e^{-t\Sigma}] = \exp\left\{ -\int_S (1-e^{-tg(x)})\mu(\dx) \right\}\,.
\end{equation}

\subsection{Completely random measures}~\label{section:CRM}
Let $(\Omega, \mathcal{F}, \mathbb{P})$ be some probability space. Let $(M(S), \mathcal{B})$ be the space of all $\sigma$-finite measures on $(S, \mathcal{S})$ supplied with an appropriate $\sigma$-algebra. A completely random measure (CRM) $\Lambda$ on $(S, \mathcal{S})$, is a measurable mapping from $\Omega$ to $M(S)$ such that
\begin{enumerate}
\item $\mathbb{P}\{\Lambda(\emptyset)=0 \} = 1$,
\item For any disjoint countable collection of sets $A_1, A_2,\dotsc,$ the random variables $\Lambda(A_i), i=1,2,\dots$ are independent, and $\Lambda(\cup A_i) = \sum_i \Lambda(A_i)$, holds almost surely. (the independent increments property)
\end{enumerate}

An important characterization of CRMs in terms of Poisson processes is as follows~\cite{kingman1967completely}. For any CRM $\Lambda$ on $(S, \mathcal{S})$ without any fixed atoms or deterministic component, there exists a Poisson process $N$ on $(\mathbb{R}^+\times S , \mathcal{B}_{\mathbb{R}^+} \otimes \mathcal{S})$, such that $\Lambda(\dx) = \int_{\mathbb{R}^+}zN(\dz, \dx)$. Using Campbell's proposition, the Laplace transform of $\Lambda(A)$ for a measurable set $A$, is given by the following formula:
\begin{equation}
\mathbb{E}[e^{-t\Lambda(A)}] = \exp\left(-\int_{\mathbb{R}^+\times A} (1-e^{-tz}) \nu(\mathrm{d}z, \mathrm{d}x)\right),\, t\ge 0\:,
\end{equation} 
where $\nu$ denotes the mean measure of the underlying Poisson process $N$. $\nu$ is also referred to as the Poisson intensity measure of $\Lambda$. If $\nu(\mathrm{d}z, \mathrm{d}x) = \rho(\mathrm{d}z)\mu(\mathrm{d}x)$, for a $\sigma$-finite measure $\mu$ on ${S}$, and a $\sigma$-finite measure $\rho$ on $\mathbb{R}^+$ that satisfies $\int_{\mathbb{R}^+}(1-e^{-tz})\rho(\dz) < \infty$, then $\Lambda(.)$ is known as homogenous CRM. In sequel, we assume $\mu(.)$ to be finite. Moreover, unless specified, whenever we refer to CRM, it means a homogeneous completely random measure without any fixed atoms or deterministic component.

Let $N$ be the Poisson process of the CRM $\Lambda$, that is, $\Lambda(\dx) = \int_{\mathbb{R}^+} s N(\dz, \dx)$. If $\Pi$ is the random collection of points corresponding to $N$, then $\Lambda$ can equivalently be written as $\Lambda = \sum_{(z, x)\in \Pi} z \delta_{x}$. $\{z:(z,x)\in \Pi\}$ constitute the weights of the CRM $\Lambda$. By the mapping proposition for Poisson processes, they form a Poisson process with mean measure $\mu^*(\dz) = \mu \circ f^{-1}(\dz)$, where $f(x,y) = x$ is the projection map on $\mathbb{R}^+$. Hence, the weights of $\Lambda$ form a Poisson process on $\mathbb{R}^+$ with mean measure $\mu^*(\dz) = \nu(\dz, S) = \rho(\dz)\mu(S)$. We formally state this result below.

\begin{lemma}~\label{CRMweights}
The weights of a homogenous CRM  with no atoms or deterministic component, whose Poisson intensity measure $\nu(\dz, \dx) = \rho(\dz)\mu(\dx)$ form a Poisson process with mean measure $\rho(\dz)\mu(S)$.
\end{lemma}

\textbf{Note 1:} A completely random measure without any fixed atoms or deterministic component is a purely-atomic random measure.

\textbf{Note 2:} Every such homogeneous CRM $\Lambda$ on $(S, \mathcal{S})$ has an underlying Poisson process $N$ on $(\mathbb{R}^+\times S, \mathcal{B}_{\mathbb{R}^+} \otimes \mathcal{S} )$, such that
\begin{equation}
\Lambda(\dx) = \int_{\mathbb{R}^+} z N(\dz, \dx)
\end{equation}
almost surely.

\section{The proposed model}~\label{model}
Let $X_1,\dotsc,X_n$ be $n$ observed samples, for instance, $n$ documents. We assume that each sample $X_i$ is generated as follows:
\begin{itemize}
\item The base measure $\Phi$ is CRM$(\rho, \mu)$, where $\rho$ and $\mu$ are $\sigma$-finite and finite (non-atomic) measures on $(S, \mathcal{S})$ respectively.
\item Object specific measures $\Lambda_i, 1\le i\le n$ are CRM$(\bar{\rho}, \Phi)$, where $\bar{\rho}$ is another $\sigma$-finite non-atomic measure on $(S, \mathcal{S})$.
\item The latent feature set $N_i$ for each object $X_i$ is a Poisson process with mean measure $\Lambda_i$. 
\item Finally, the visible features $X_i$ are sampled from $N_i$.
\end{itemize}

\textbf{Note:} For topic modelling, $S$ corresponds to the space of all probability measures on the words in the dictionary, also known as topics. Hence, when we sample $\Phi$, we sample a subset of topics, along with the weights for those topics. This follows from the discreteness of $\Phi$. Sampling object-specific random measures $\Lambda_i$ corresponds to sampling the document specific weights for  all the topics in $\Phi$. Sampling the latent features $N_i $ then corresponds to selecting a subset of topics from $\Lambda_i$ based on the corresponding document-specific weights. Since, all the $\Lambda_i's$ have access to the same set of topics, this leads to sharing of topics among $N_i$s. Finally, the words in $X_i$ is sampled from the corresponding topic in $N_i$ using categorical distribution.

Our aim is to infer the latent features $N_i, 1\le i \le n$ from $X_i, 1\le i \le n$. By Bayes' rule
\begin{align*}
P&(N_1 , \dotsc, N_n|X_1, \dotsc, X_n) \propto \\ &P(X_1, \dotsc, X_n|N_1, \dotsc, N_n) P(N_1, \dotsc, N_n) \\
&=\Pi_{i=1}^n P(X_i | N_i) P(N_1, \dotsc, N_n)
\end{align*} 
The conditional distribution of $X_i$ given $N_i$ are often very simple to compute, for instance, in the case of topic modelling, it is simply the product of categorical distributions. Hence, all we need to compute is the prior distribution of the latent features $N_1, \dotsc, N_n$. This can be obtained by marginalizing out the base and object-specific random measures $\Phi$ and $\Lambda_i, 1\le i \le n$. This is what we wish to achieve in the next few sections.

We will address the problem of marginalizing out the base and object-specific random measures in two steps. Firstly, in section~\ref{object_specific}, we will derive results for the case when the base measure is held fixed and the object-specific random measure is marginalized out. Next, in section~\ref{base}, we will derive results for the case, when the base random measure $\Phi$ is also marginalized out. All the proofs are provided in the appendix.

\subsection{Marginalizing out the object specific measure}~\label{object_specific}
Let $\phi$ be a realization of the base random measure $\Phi$. Let $\Lambda_i, 1 \le i \le n$, be independent CRM$(\bar{\rho}, \phi)$. It is straightforward to see that if $\phi$ is a finite measure $\Lambda_i$s will almost-surely be finite. Because of the independence among $\Lambda_i$s, we can focus on marginalizing out a single object-specific random measure, say $\Lambda$. Although, in our original formulation, only $1$ object is sampled from its object-specific random measure, we will present results for the case when $n$ objects, $N_1,\dotsc, N_n$ are sampled from the object specific random measure. This extended result will be needed in the next section when marginalizing the base measure.

There are several ways to instantiate the random measure $\Lambda$. For instance, one can use the fact that since the underlying base measure $\phi$ is purely-atomic, the support of CRM$(\bar{\rho}, \phi)$ will be restricted to only those measures whose support is a subset of the support of $\phi$. In particular, if $\phi = \sum_{j=1}^\infty \beta_j \delta_{x_j}$, then $\Lambda$ will be of the form $ \sum_{j=1}^\infty L_j \delta_{x_j}$, where $L_j$ are independent random variables. The independence of $L_j$s follows from the complete randomness of the measure.

However, we found that this approach doesn't lead us far. Hence, we derive the marginal distribution of the Poisson processes $N_1, \dotsc, N_n$ in proposition~\ref{CRMdistinct} and~\ref{CRMPoisson}, by first assuming $\phi$ to be a continuous measure and then generalizing it to the case where $\phi$ is any finite measure. 

In the sequel, $\psi(t) = \int_{\mathbb{R}^+}(1-e^{-tz})\bar{\rho}(\dz)$, and $\psi^{(k)}$ is the $k^{th}$ derivative of $\psi$.

\begin{proposition}~\label{CRMdistinct}
Let $\Lambda$ be a CRM on $(S, \mathcal{S})$ with Poisson intensity measure $\rho(\dz)\mu(\dx)$, where both $\mu(.)$ and $\rho(.)$ are non-atomic. Let $N_1,\dotsc, N_n$ be n independent Poisson process with random mean measure $\Lambda$, and $M$ be the distinct points of $N_i,\, 1\le i\le n$. Then, 
$M$ is a Poisson process with mean measure $\mathbb{E}[M(\dx)] = {\mu(\dx)}\int_{\mathbb{R}^+}(1-e^{-nz})\rho(\dz)$.
\end{proposition}

The above proposition provides the distribution of distinct points of the $n$ point processes, $N_1,\dotsc, N_n$. In order to complete the description of the distribution of $N_1,\dotsc, N_n$, we also need to specify the joint distribution of the counts of each distinct feature in each $N_i$. This distribution is referred to as CRM-Poisson distribution in the rest of the paper. Let $M(S) = k$ and $m_{ij}$ be the count of the $j^{th}$ distinct feature in the $i^{th}$ object. Furthermore, let $[m_{\LargerCdot j}]$ be the count of the $j^{th}$ distinct feature for each object and $[m_{ij}]_{1\le i\le n, 1 \le j \le k}$ be the set of count vectors for the each latent feature.

\begin{proposition}~\label{CRMPoisson}
The joint distribution of the set of count vectors for the each latent feature $[m_{ij}]_{(n,k)}$ is given by
\begin{equation} \label{ESPF}
P([m_{ij}]_{(n,k)}) = {(-1)^{m_{\LargerCdot \LargerCdot}-k}}\frac{\theta^k e^{-\theta\psi(n)}}{\prod_{i=1}^n( m_{i\LargerCdot})!} \prod_{j=1}^k \psi^{(m_{\LargerCdot j})}(n)\,,
\end{equation}
where $m_{i\LargerCdot} = \sum_{j=1}^k m_{ij}$, $m_{\LargerCdot j} = \sum_{i=1}^n m_{ij}$, $m_{\LargerCdot \LargerCdot} = \sum_{i=1}^n \sum_{j=1}^k m_{ij}$, $\theta = \mu(S)$ and $\psi(t) =  \int_{\mathbb{R}^+} (1-e^{-tz})\rho(\dz) $ is the Laplace exponent of $\Lambda$,  and $\psi^{(l)}(t)$ is the $l^{th}$ derivative of $\psi(t)$. This distribution will be referred to as \textbf{CRM-Poisson}$(\mu(S), \rho, n)$.

\end{proposition}

\begin{corollary}~\label{CRMPoissonCon}
Conditioned on $M(S)=k$, the set of count vectors for the each latent feature $ [m_{ij}]_{(n, k) }$ is distributed as 
\begin{align} 
&P([m_{ij}]_{(n,k)}| M(S)=k) \label{eq:CRMPoissonCon} \\ &= \frac{\theta^k(-1)^{m_{\LargerCdot \LargerCdot}-k}k!}{\prod_{i=1}^n(m_{i\LargerCdot})!} \prod_{j=1}^k \frac{\psi^{(m_{i\LargerCdot})}(n)}{\psi(n)} \notag
\end{align} 
\end{corollary}

Note that both $\psi^{(k)}$ and $\psi$ contain a multiple involving $\mu(S)$, which cancels out when they are divided in~\eqref{eq:CRMPoissonCon}. Hence, conditioned on the number of points in the Poisson process $M$, the distribution of the set of counts for each latent feature $[m_{ij}]_{(n, k) }$ does not depend on the measure $\mu$. In sequel, this distribution will be referred to as \textbf{conditional CRM-Poisson}$(\rho, n, k)$ or \textbf{CCRM-Poisson}$(\rho, n, k)$.

\textbf{Example 1: The Gamma-Poisson process}\\
The Poisson-intensity measure of gamma process is given by $\rho(\dz) = e^{-z}z^{-1}\dz$. The corresponding Laplace exponent is $\psi(t) = \ln(1+t)$. Replacing it in equation~\eqref{ESPF}, we get
\begin{equation}
P([m_{ij}]_{(n,k)}) = \frac{\theta^k \prod_{j=1}^k \Gamma(m_{\LargerCdot j})}{\prod_{i=1}^n m_{i\LargerCdot}!(1+n)^{m_{\LargerCdot \LargerCdot}+\theta}} 
\end{equation}

Next, we generalize these results for the case when $\phi$ is an atomic measure.
\begin{proposition}~\label{CRMnonAtomic}
Let $\Lambda$ be a completely random measure with Poisson intensity measure $\nu(\dz, \dx) = \phi(\dx)\bar{\rho}(\dz)$, where $\bar{\rho}$ is non-atomic. Let $N$ be a Poisson process with mean measure $\Lambda$. Then, $N$ can be obtained by sampling a Poisson process with mean measure $\phi(\dx)\psi(1)$, say $M$, and then sampling the count of each feature in $M$ using the conditional CRM-Poisson distribution.
\end{proposition}
\textbf{Note:} The points in $M$ won't be distinct anymore, since the underlying mean measure is non-atomic.

\subsection{Marginalizing out the base measure}~\label{base}
The previous section derived the marginal distribution of the Poisson processes, for a fixed realization $\phi$ of the base random measure $\Phi$. In this section, we want to marginalize the CRM $\Phi$ as well. Marginalizing $\Phi$ does away with the independence among the latent features $N_i$s, hence, we need to model the joint distribution of $N_1, \dotsc, N_n$.

The model under study is
\begin{align}~\label{CRMPP1}
\begin{split}
\Phi &\sim \text{CRM}(\rho, \mu) \:,\\
\Lambda_i|\Phi &\sim \text{CRM}(\rho^{\prime}, \Phi),\, 1\le i \le n \:, \\
N_i|\Lambda_i &\sim \text{Poisson Process}(\Lambda_i),\, 1\le i \le n \:.
\end{split}
\end{align}
We use Proposition~\ref{CRMnonAtomic} to marginalize out $\Lambda_i$ from the above description. Thus $N_i$ can equivalently be obtained by sampling a Poisson processes with mean measure $\Phi(\dx)\int_{\mathbb{R}^+}(1-e^{-z})\rho(\dz)$, and then sampling the count of each feature in $M_i$ for each point process $N_i$ using Corollary~\ref{CRMPoissonCon}. In particular, let $M_i$ be the corresponding Poisson process, and $m_{ij}$ be the count of the $j^{th}$ feature in $M_i$ for the point process $N_i$ and $r_{i\LargerCdot} = M_i(S)$. The reason for the symbol $r_{i\LargerCdot}$ will become clear, when we have a picture of the entire generative model. Let $[m_{ij}]_{\LargerCdot, r_{i\SmallerCdot}}$ be the set of counts of the latent features for the $i^{th}$ individual. The distribution of the set of counts $[m_{ij}]_{\LargerCdot, r_{i\SmallerCdot}}$ conditioned on $M_i(S)$ does not depend on $\Phi$. Hence, an alternative description of the $N_i$ via $M_i$ and $m_{ij},\,1\le j\le r_{i\LargerCdot}$ is as follows:
\begin{align}~\label{CRMPP2}
\begin{split}
&{M_i}|\Phi \sim \text{Poisson Process}\left(\Phi(.)\int_{\mathbb{R}^+}(1-e^{-z})\bar{\rho}(\dz)\right)\:, \\
&[m_{ij}]_{(\LargerCdot, r_{i\SmallerCdot})}| \{M_i(S)=r_{i\LargerCdot} \} \sim \text{CCRM-Poisson}(\bar{\rho},1,r_{i\LargerCdot})\, \\
& N_i = \sum_{i=1}^{r_{i\SmallerCdot}} m_{ij} \delta_{M_{ij}}\,,
\end{split}
\end{align}
where $M_{ij}$ are the points in the point process $M_i$.

$M_i,\, 1\le i \le n$ are independent Poisson processes, whose mean measure is a scaled CRM, and hence, also a CRM. Hence, we are again in the domain of CRM-Poisson models. Let $\bar{\psi}(1) = \int_{\mathbb{R}^+}(1-e^{-z})\bar{\rho}(\dz)$. If we define $\Phi'(\dx) = \bar{\psi}(1)\Phi(\dx)$, then 
\begin{align*}
&\mathbb{E}[e^{-t\Phi'(A)}] = \mathbb{E}[e^{-t\bar{\psi}(1)\Phi(A)}] \\&= \exp\left\{-\mu(A)\int_{\mathbb{R}^+} (1-e^{-t\bar{\psi}(1) z}) \rho(\mathrm{d}z)\right\} \\
& = \exp\left\{-\mu(A)\int_{\mathbb{R}^+} (1-e^{-tz'}) \rho(\mathrm{d}(z'/\bar{\psi}(1)))\right\}
\end{align*}
Hence, the Poisson intensity measure of the scaled CRM $\Phi'$ is given by $\rho(\mathrm{d}(z/\bar{\psi}(1)))\mu(\dx)$. Applying Proposition~\ref{CRMnonAtomic} to marginalize out $\Phi$, we get that $M_i$'s can be obtained by sampling a Poisson process $R$ with mean measure 
\begin{align*}
&\mathbb{E}[R(\dx)] = \mu(\dx)\int_{\mathbb{R}^+}(1-e^{-nz'})\rho(\mathrm{d}(z'/\bar{\psi}(1))) \\&= \mu(\dx)\int_{\mathbb{R}^+}(1-e^{-\bar{\psi}(1) nz})\rho(\mathrm{d}z)\,.
\end{align*}
The count of each feature in $R$ for each point process $M_i$ can then be obtained by using Corollary~\ref{CRMPoissonCon}. In particular, let $r_{ik}$ be the count of the $k^{th}$ point in $R$ for the point process $M_i$ and $p = R(S)$.

A complete generative model for generating the point processes $N_i,\, 1\le i\le n$ is as follows:
\begin{align}
&{R} \sim \text{Poisson Process}\left(\mu(.)\int_{\mathbb{R}^+}(1-e^{-\bar{\psi}(1) nz})\rho(\dz)\right),~\notag \\
&[r_{ik}]_{(n, p) }|\left\{ R(S)=p\right\} \sim \text{CCRM-Poisson}(\rho,\bar{\psi}(1) n, p)\,~\label{CRMPP3}\\
&M_i = \sum_{k=1}^p r_{ik}\delta_{R_k}~\notag\\
&[m_{ij}]_{(\LargerCdot, r_{i\SmallerCdot})}| \{M_i(S)=r_{i\LargerCdot}\} \sim \text{CCRM-Poisson}(\bar{\rho},1,r_{i\LargerCdot})  ~\notag\\
&N_i = \sum_{j=1}^{r_{i\SmallerCdot}} m_{ij} \delta_{M_{ij}}\,,\notag
\end{align}
Since $R$ is again a Poisson process, it is straightforward to extend this hierarchy further by sampling $\mu(.)$ again from a CRM.

\section{Implementation via Gibbs sampling}
Section~\ref{model} provided an approach for sampling a Poisson process, when sampled from a hierarchical CRM, without having to instantiate the infinitely many atoms of the base or object-specific CRM. However, it is not clear how the above derivations can be used for determining the latent features $N_1, \dotsc, N_n$ for the objects $X_1, \dotsc, X_n$, which is the aim of this work.

In this section, we provide a Gibbs sampling approach for sampling the latent features from its prior distribution that is $P(N_1, \dotsc, N_n)$. In order to sample from the posterior, one simply needs to multiply the equations in this section with the likelihood of the latent feature. In order to be able to perform MCMC sampling in hierarchical CRM-Poisson models, we need to marginalize out $R(S)$ and $M_i(S)$ from distributions of $[r_{ik}]_{(n, p) }$ and $[m_{ij}]_{(\LargerCdot, r_{i\SmallerCdot})}$ respectively. By marginalizing out the Poisson distributed random variable $R(S)$ from~\eqref{CRMPP3}, we get that 
\begin{equation*}
[r_{ik}]_{(n, p) } \sim \text{CRM-Poisson}(\mu(S),\rho,\bar{\psi}(1) n)\,.
\end{equation*}
The marginal distribution of the set of counts of each latent feature for the $i^{th}$ individual $[m_{ij}]_{(\LargerCdot, r_{i\SmallerCdot})}$ (where $r_{i\SmallerCdot}$ is also random) is given by the following lemma.
\begin{lemma}~\label{tables}
Let
$$h(u) = \mathbb{E}[e^{-u\psi(S)}] = \exp\left\{-\mu(S)\int_{\mathbb{R}^+}(1-e^{-uz})\rho(\dz)\right\}.$$ Furthermore, if we let $$\psi(u) = \int_{\mathbb{R}^+}(1-e^{-uz})\rho(\dz)$$ $$\bar{\psi}(u) = \int_{\mathbb{R}^+}(1-e^{-uz})\bar{\rho}(\dz)\, ,$$ 
then, $[m_{ij}]_{(\LargerCdot, r_{i\SmallerCdot})}$ is marginally distributed as 
\begin{equation}~\label{ESPF2}
P([m_{ij}]_{(\LargerCdot, r_{i\SmallerCdot})}) =  (-1)^{m_{i\LargerCdot}} h^{(r_i\SmallerCdot)}\left(\bar{\psi}(1)\right) \frac{\prod_{j=1}^{r_{i\SmallerCdot}} \bar{\psi}^{(m_{ij})}(1)}{ m_{i\SmallerCdot}!}\,,
\end{equation}
\end{lemma}

In the case of topic-modelling, the number of latent features, $\# N_i$ is equal to the number of observed features $\# X_i$. Hence, let $X_{il}$ be the $l^{th}$ observed feature associated with the $i^{th}$ object and $N_{il}$ be the corresponding latent feature. Here, we discuss the MCMC  approach for sampling from the prior distribution of $N_{il}, \, 1\le l \le m_{i\LargerCdot}$.

As discussed in~\cite{neal2000markov}, it is more efficient to sample the index of the latent feature, rather than the latent feature itself. Hence, let $T_{il}$ be the index of the point in $M_i$ associated with $N_{il}$, and $D_{ij}$ be the index of the point in $R$ associated with $M_{ij}$. In an analny with the Chinese restaurant franchise model~\cite{teh2006hierarchical}, one can think of $T_{il}$ to be the index of the table assigned to the $l^{th}$ customer in the $i^{th}$ restaurant, and $D_{ij}$ to be the index of the dish associated with the $j^{th}$ table in $i^{th}$ restaurant. Moreover, $m_{ij}$ refers to the number of customers sitting on the $j^{th}$ table in $i^{th}$ restaurant, and $r_{ik}$ refers to the number of tables in the $i^{th}$ restaurant with the $k^{th}$ dish. Hence $r_{i\LargerCdot} = \sum_{k=1}^{p}r_{ik}$ is the number of tables in the $i^{th}$ restaurant.

The distribution of the number of customers per table in the $i^{th}$ restaurant, $[m_{ij}]_{(\LargerCdot, r_{i\SmallerCdot})}$ follows from Lemma~\ref{tables}. Hence, in order to sample the table of $l^{th}$ customer, $T_{il}$, given the indices of the tables of all the other customers in $i^{th}$ restaurant, we treat it as the table corresponding to the last customer of the $i^{th}$ restaurant. Let $m_{i'j}^{-(il)}$ be the number of customers sitting on the $j^{th}$ table in the ${i'}^{th}$ restaurant, excluding the $l^{th}$ customer. The probability that the $l^{th}$ customer in the $i^{th}$ restaurant occupies the $j^{th}$ table is proportional to $P({m}_{ij'}^{-(il)}+ 1_{j' = j}, 1\le j' \le r_{i\LargerCdot})$ as given in~\eqref{ESPF2}. We divide the expression by $P({m}_{ij'}^{-(il)}, 1\le j' \le r_{i\LargerCdot})$ to get a simpler form for the unnormalized probability distribution.
Hence, the probability of assigning an existing table with index $j$ is given by
\begin{equation}
P(T_{il} = j| \mathbf{T}^{-(il)}) \propto -\frac{  \bar{\psi}^{ (m_{i j}^{-(il)}+1)}(1)}{\bar{\psi}^{ (m_{i j}^{-(il)})}(1)}\,,\label{eq:tableOld}
\end{equation}
and the probability of sampling a new table for the customer is given by 
\begin{equation}
P(T_{il} = r_{i\LargerCdot}+1| \mathbf{T}^{-(il)}) \propto =-\frac{h^{(r_{i\LargerCdot}+1)}(\bar{\psi}(1))}{h^{(r_{i\LargerCdot})}(\bar{\psi}(1))} \bar{\psi}^{ (1)}(1)\,,\label{eq:tableNew}
\end{equation}
where $\bar{\psi}(t)= \int_{\mathbb{R}^+}(1-e^{-tz})\bar{\rho}(\dz)$ and $h^{(k)}$ is the $k^{th}$ derivative of $h$.

Moreover, whenever a new table is sampled for a customer, a dish is sampled for the table from the distribution on tables per dish. By the discussion in the beginning of this section, the number of tables per dish $[r_{ik}]_{(n, p) }$ follow a {CRM-Poisson}$(\mu(S),\rho,\bar{\psi}(1) n)$ distribution. Hence, in order to sample the dish at $j^{th}$ table, $D_{ij}$, given the indices of the dishes at all the other tables, we treat it as the dish corresponding to the last table of the last restaurant. Let $r_{\LargerCdot k}^{-(ij)}$ be the total number of tables served with the $k^{th}$ dish, excluding the $j^{th}$ table of $i^{th}$ restaurant. The probability that the $k^{th}$ dish  is served at the $j^{th}$ table in the $i^{th}$ restaurant is proportional to $P({r}_{i'k'}^{-(ij)}+ 1_{i' = i', k'=k}, 1 \le i' \le n, 1\le k' \le p)$ as given in~\eqref{ESPF2}. We divide the expression by $P({r}_{i'k'}^{-(ij)}, 1 \le i' \le n, 1\le k' \le p)$ to get a simpler form for the unnormalized probability distribution. Hence, the probability of serving an existing dish with index $k$ is given by
\begin{equation}
P(D_{ij} = k| \mathbf{D}^{-(ij)}) \propto - \frac{  \psi^{ (r_{\LargerCdot k}^{-(ij)}+1)}(\bar{\psi}(1) n)}{\psi^{ (r_{\LargerCdot k}^{-(ij)})}(\bar{\psi}(1) n)}\,,\label{eq:dishOld}
\end{equation}
and the probability of sampling a new dish for the table is given by 
\begin{equation}
P(D_{ij} = p+1| \mathbf{D}^{-(ij)}) \propto  \theta \psi^{ (1)}(\bar{\psi}(1) n)\,,\label{eq:dishNew}
\end{equation}
where $\psi(t)= \int_{\mathbb{R}^+}(1-e^{-tz})\rho(\dz)$ and $\theta = \mu(S)$.

Hence, a complete description of one iteration of MCMC sampling, from the prior distribution, in hierarchical CRM-Poisson models is as follows:
\begin{enumerate}
\item For each customer in each restaurant, sample his table index conditioned on the indices of table of other customers, according to equations~\eqref{eq:tableOld} and~\eqref{eq:tableNew}.

\item If the table selected is a new table, sample the index of dish corresponding to that table from equations~\eqref{eq:dishOld} and~\eqref{eq:dishNew}.

 \item Sample the index of dish for each table, conditioned on the indices of dishes at the other tables, according to equations~\eqref{eq:dishOld} and~\eqref{eq:dishNew}.
 \end{enumerate}

\textbf{Example 2: The Gamma-Gamma-Poisson process}\\
We compute the dish and table sampling probabilities for the Gamma-Gamma-Poisson process using the above equations. The Poisson intensity measure for both the base and object specific measures $\Phi$ and $\Lambda_i,1\le i\le n$ is $z^{-1}e^{-z}\dz$. The corresponding Laplace exponent is given by $\psi(t) = \bar{\psi}(t) = \ln(1+t)$. Moreover, let the mean measure for the base measure $\Phi$ be $\mu(.)$ and $\mu(S) = \theta$. Then, $h(u) = \mathbb{E}e^{-u\Phi(S)} = \frac{1}{(1+u)^\theta}$. The corresponding derivatives are given by
\begin{align}
&\psi^{(k)} = \bar{\psi}^{(k)}(t) = \frac{(-1)^{k-1}{\Gamma(k)}}{(1+t)^k} \\
&h^{(k)}(u) = \frac{(-1)^k \Gamma(k+\theta)}{(1+u)^{k+\theta}\Gamma(\theta)}
\end{align}
The corresponding dish sampling probabilities are given by
\begin{align}
P(D_{ij} = k| \mathbf{D}^{-(ij)}) \propto \frac{r_{\LargerCdot k}^{-(ij)}}{1+n\ln2} \\
P(D_{ij} = p+1| \mathbf{D}^{-(ij)}) \propto \frac{\theta}{1+n\ln2}
\end{align}
for an existing and new dish respectively. Normalizing these probabilities, we get
\begin{align}
P(D_{ij} = k| \mathbf{D}^{-(ij)}) = \frac{r_{\LargerCdot k}^{-(ij)}}{\sum_{k=1}^p r_{\LargerCdot k} + \theta} \\
P(D_{ij} = p+1| \mathbf{D}^{-(ij)}) = \frac{\theta}{\sum_{k=1}^p r_{\LargerCdot k} + \theta}
\end{align}

The table sampling probabilities are given by
\begin{align}
P(T_{il} = j| \mathbf{T}^{-(il)}) \propto \frac{m_{ij}^{-(il)}}{1+\ln(2)} \\
P(T_{il} = r_{i\LargerCdot}+1| \mathbf{T}^{-(il)}) \propto \frac{\theta + r_{i\LargerCdot}}{(1+\ln(2))^2}
\end{align}
for an existing and new table respectively. Normalizing these probabilities, we get
\begin{align}
P(T_{il} = j| \mathbf{T}^{-(il)}) = \frac{m_{ij}^{-(il)}}{ \sum_{j=1}^{r_{i\LargerCdot}} m_{ij}^{-(il)} + \frac{\theta+r_{i\LargerCdot}}{1+\ln(2)}} \\
P(T_{il} = r_{i\LargerCdot}+1| \mathbf{T}^{-(il)}) = \frac{(\theta + r_{i\LargerCdot})/(1+\ln(2))}{ \sum_{j=1}^{r_{i\LargerCdot}} m_{ij}^{-(il)} + \frac{\theta+r_{i\LargerCdot}}{1+\ln(2)}}
\end{align}

\textbf{Example 3: The Gamma-Generalized Gamma-Poisson process}\\
In this scenario, the base random measure has $\rho(\dz) = e^{-z}z^{-1}\dz$, whereas the object specific measure has $\bar{\rho}(\dz) = e^{-z} z^{-d-1} \dz$, where $0<d<1$ is known as the discount parameter. The corresponding Laplace exponents are given by $\psi(t) = \ln(1+t)$ and $\bar{\psi}(t) = \frac{(1+t)^d-1}{d}$ respectively. The derivative of $\bar{\psi}$ is given by
\begin{align}
&\bar{\psi}^{(k)}(t) = \frac{(-1)^{k-1}{\Gamma(k-d)}}{(1+t)^{k-d} \Gamma(1-d)} \\
\end{align}
Other derivatives remain same as in the previous example. Moreover, the dish sampling probabilities remain same. The table sampling probabilities are given by
\begin{align}
P(T_{il} = j| \mathbf{T}^{-(il)}) = \frac{m_{ij}^{-(il)} -d}{ \sum_{j=1}^{r_{i\LargerCdot}} (m_{ij}^{-(il)}-d) + \frac{\theta+r_{i\LargerCdot}}{1+\ln_d(2)}} \\
P(T_{il} = r_{i\LargerCdot}+1| \mathbf{T}^{-(il)}) = \frac{(\theta + r_{i\LargerCdot})/(1+\ln_d(2))}{ \sum_{j=1}^{r_{i\LargerCdot}} (m_{ij}^{-(il)}-d) + \frac{\theta+r_{i\LargerCdot}}{1+\ln_d(2)}}
\end{align}
where $\ln_d(2) = \frac{2^d-1}{d}$.
\section{Experimental results}
We use hierarchical CRM-Poisson models for learning topics from the NIPS corpus~\footnote{The dataset can be downloaded from \url{http://psiexp.ss.uci.edu/research/programs_data/toolbox.htm}}.

\subsection{Evaluation} For evaluating the different models, we divide each document into a training section and a test section by independently sampling a boolean random variable for each word. The probability of sending the word to the training section is varied from $0.3$ to $0.7$. We run $2000$ iterations of Gibbs sampling. The first $500$ iterations are discarded, and every sample in every $5$ iterations afterwards is used to update the document-specific distribution on topics and the topic specific distribution on words. In particular, let $W$ be the number of words, $K$ be the number of topics, $(\beta_{dk})_{1\le k\le K}$ be the document specific distribution on topics for the document $d$, and $(\tau_{kw})_{1\le w\le W}$ be the topic specific distribution on words for the $k^{th}$ topic. Then, the probability of observing a word $w$ in document $d$ is given by $\sum_{k=1}^K \beta_{dk}\tau_{kw}$. For the evaluation metric, we use perplexity, which is simply the inverse of the geometric mean of the probability of all the words in the test set.

\subsection{Varying the Common CRM}
In our experiments, we fix the object specific random measure $\Lambda_i$ in~\eqref{CRMPP1} to be the gamma process, with $\bar{\rho}(\dz) = e^{-z}z^{-1}\dz$. For the base CRM $\Phi$, we consider two specific choices of random measures.
\begin{itemize}
\item 
\textbf{Generalized gamma process (GGP):} The Poisson intensity measure of $\Phi$ is given by $\nu(\dz, \dx) = \rho(\dz)\mu(\dx)$, where $\rho(\dz) = \frac{\theta}{\Gamma(1-d)} e^{-z}z^{-d-1}\dz,\, 0\le d<1, \theta > 0$ and $\mu(S) = 1$. The corresponding Laplace exponent is given by $\theta({(1+{t})^d-1})/{d}$.

\item 
\textbf{Sum of Generalized gamma processes (SGGP):}
The Poisson intensity measure of the CRM is given by $\nu(\dz, \dx) = \rho(\dz)\mu(\dx)$, where 
\begin{equation}
\rho(\dz) = \sum_{q=1}^m \frac{\theta_q}{\Gamma(1-d_q)} e^{-z}z^{-d_q-1}\dz \,
\end{equation}
and $\mu(S) = 1$. The corresponding Laplace exponent is given by 
\begin{equation}
\psi(t) = \left(\sum_{q=1}^m\theta_q \frac{(1+t)^{d_q} - 1}{d_q}\right)\,.
\end{equation}
 \end{itemize}

For the case of GGP, the value of the discount parameter $d$ is chosen from the set $\{0,\,.1,\,.2,\,.3,\,.4\}$. Furthermore, a gamma prior with rate parameter $2$ and shape parameter $4$ is defined on $\theta$.

\textbf{Note:} The generalized gamma process with discount parameter $0$ corresponds to the Gamma process. Using a gamma process prior for the base and object-specific CRM corresponds exactly to the hierarchical Dirichlet process with a gamma prior on the concentration parameter of the object specific Dirichlet process. We did not add comparison results with HDP separately, because the same perplexity is obtained in both the models.

\begin{figure}[h]
\begin{center}
\includegraphics[width=.45\textwidth]{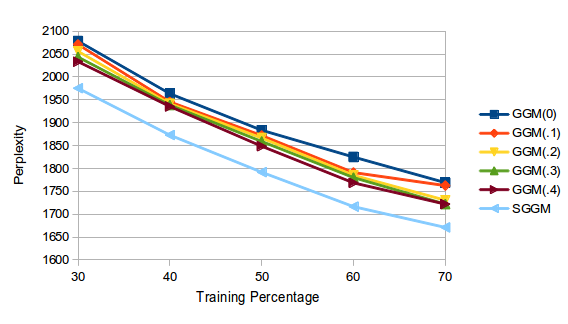}
\end{center}
\caption{Variation of perplexity with training percentage for various hierarchical CRM-Poisson models}
\label{fig:perplexity}
\end{figure}

For the case of SGGP, we consider $m=5$, and $d_1 = 0,\, d_2=.1\dotsc, d_5=.4$. Furthermore, independent gamma priors with rate parameter $2$ and shape parameter $4$ are defined for each $\theta_q,\, 1\le q \le 5$. The posterior of each parameter $\theta_q$ is sampled via uniform sampling.
We use equations~\eqref{eq:tableOld}-\eqref{eq:dishNew} to compute the dish sampling and table sampling probabilities. The probability of sampling an existing dish is given by
\begin{align*}
& P(D_{ij} = k| \mathbf{D}^{-(ij)}) \\
& \propto \frac{\sum_{q=1}^m \theta_q \frac{\Gamma(r_{\LargerCdot k}^{-(ij)}+1-d_q)}{\Gamma(1-d_q)}(1+\bar{\psi}(1) n)^{d_q}}{\sum_{q=1}^m \theta_q \frac{\Gamma(r_{\LargerCdot k}^{-(ij)}-d_q)}{\Gamma(1-d_q)}(1+\bar{\psi}(1) n)^{d_q}}\,,
\end{align*}
where $\bar{\psi}(1) = \int_{\mathbb{R}^+}(1-e^{-z})\bar{\rho}(\dz) = \ln(2)$. Similarly, the probability of a new dish is given by
\begin{equation*}
P(D_{ij} = p+1| \mathbf{D}^{-(ij)}) \propto   \sum_{q=1}^m \theta_q (1+\bar{\psi}(1) n)^{d_q}\,.
\end{equation*}
The table-sampling probabilities can be computed similarly. We approximated the Laplace transform of $\Phi(S)$~($h$ in \eqref{eq:tableNew}), by a weighted sum of exponential functions to simplify the computation of its derivatives. The perplexity for the hierarchical CRM-Poisson models as a function of training percentage is plotted in Figure~\ref{fig:perplexity}. Note that Figure~\ref{fig:perplexity} doesn't necessarily imply that SGGM-based models will always outperform GGM based models as the results have been obtained by defining a specific gamma prior for each hyperparameter, as mentioned above.

\section{Conclusion}
For years, hierarchical Dirichlet processes have been the standard tool for nonparametric topic modelling, since collapsed inference in HDP can be performed using the Chinese restaurant franchise scheme. In this paper, our aim was to show that collapsed Gibbs sampling can be extended to a much larger set of hierarchical random measures using the same Chinese restaurant franchise scheme, thereby opening doors for further research into the efficacy of various hierarchical priors. We hope that this will encourage a better understanding of applicability of various hierarchical CRM priors. Furthermore, the results of the paper can be used to prove results for hierarchical CRMs in other contexts, for instance, nonparametric hidden Markov models.

\section*{Acknowledgement}
Gaurav Pandey is supported by IBM PhD Fellowship for the academic year 2015-2016.

\bibliography{nips}
\bibliographystyle{icml2016}

\newpage
   
\newpage
\section*{Appendix}

\textbf{Proof of proposition~\ref{CRMdistinct} }
\begin{proof}
Let $N = \sum_{i=1}^n N_i$. Since, conditioned on $\Lambda$, $N_1,\dotsc, N_n$ are independent Poisson process with mean measure $\Lambda$, by the superposition proposition for Poisson processes, $N$ is a Poisson process conditioned on $\Lambda$ with mean measure $n\Lambda$. Since, $\mathbb{E}[e^{t_1N(A) + t_2 N(B)}| \Lambda] = \mathbb{E}[e^{t_1N(A)}|\Lambda(A)]\mathbb{E}[e^{t_2N(B)}|\Lambda(B)]$, and $\Lambda(A)$ and $\Lambda(B)$ are independent, hence, $N(A)$ and $N(B)$ are also independent and therefore, $N$ is a CRM. Hence, $N(\dx) = \int_{\mathbb{R}^+}s\bar{N}(\dz, \dx)$, for a Poisson process $\bar{N}$ on $S\times\mathbb{R}^+$. Moreover,
\begin{align*}
&\mathbb{E}[e^{-tN(A)}] \\
&= \mathbb{E}\left[\mathbb{E}[e^{-tN(A)}|\Lambda] \right] \\
& =  \mathbb{E}\left[ \exp\left(-n\Lambda(A)(1-e^{-t})  \right)\right] \\
& =  \exp\left(-\mu(A)\int_{\mathbb{R}^+}(1-e^{-n(1-e^{-t})z})\rho(\dz)\right) \\
& = \exp\left(-\mu(A)\int_{\mathbb{R}^+}\left(\sum_{k=0}^\infty \frac{e^{-nz}(nz)^k}{k!} \right.\right. \\
& \hspace{2cm}\left.\left. -\sum_{k=0}^\infty \frac{e^{-nz}{e^{-kt}(nz)^k}}{k!}\right)\rho(\dz)\right) \,
\end{align*}
where, we have used the fact that $1 = \sum_{k=0}^\infty\frac{e^{-nz}(nz)^k}{k!}$.
Rearranging the terms in the above equation, we get
\begin{align*}
&\mathbb{E}[e^{-tN(A)}] \\
&= \exp\left( -\mu(A)\sum_{k=1}^\infty(1-e^{-kt})\int_{\mathbb{R}^+} \frac{e^{-nz}(nz)^k}{k!}\rho(\dz) \right)\,.
\end{align*}
Hence, the Poisson intensity measure of $N$, when viewed as a CRM is given by 
\begin{equation*}
\bar{\nu}(\dk, \dx) = \mu(\dx)\int_{\mathbb{R}^+} \frac{e^{-nz}(nz)^k}{k!}\rho(\dz)
\end{equation*} 
when $k\in \{1,2,3,\dotsc\}$, and $0$ otherwise. The distinct points of $N$ can be obtained by projecting $N$ on $S$. Hence, by the mapping proposition for Poisson processes~\cite{kingman1992poisson}, the distinct points of $N$ form a Poisson process with mean measure $\mu^*(\dx) = \bar{\nu}(f^{-1}(\dx))$, where $f$ is the projection map on $S$. Hence $f^{-1}(\dx) = (\mathbb{R}^+, \dx)$, and 
\begin{align*}
&\mu^*(\dx) = \bar{\nu}(\mathbb{R}^+, \dx) \\
&= \mu(\dx)\int_{\mathbb{R}^+}\sum_{k=1}^\infty \frac{e^{-nz}(nz)^k}{k!}\rho(\dz) \\
&= \mu(\dx)\int_{\mathbb{R}^+}(1-e^{-nz})\rho(\dz)\,.
\end{align*}
Thus, the result follows.
\end{proof}

\textbf{Proof of proposition~\ref{CRMPoisson}}
\begin{proof}
The proof relies on the simple fact, that conditioned on the number of points to be sampled, the points of a Poisson process are independent~\cite{kingman1992poisson}. Thus, $n$ point processes can be sampled from a measure $\Lambda$, by first sampling the number of points in each point process from a Poisson distribution with mean $\Lambda(S)$, and then sampling the points independently. Let $\Lambda = \sum_{i=1}^n \Delta_i \delta_{X_i}$. Let $(X_{l_1}, \dotsc, X_{l_k})$ be the features discovered by the $n$ Poisson processes. Let the $i^{th}$ point process $N_i$ consist of $m_{i1}$ occurrences of $X_{l_1}$, $m_{i2}$ occurrences of $X_{l_2}$ and $m_{ik}$ occurrences of $X_{l_k}$. Then, the joint distribution of the n point processes conditioned on $\Lambda$ is given by
\begin{align*}
&\mathbb{P}(N_1,\dotsc, N_n| \Lambda) \\
&= \prod_{i=1}^n \frac{\exp(-T)T^{\sum_{j=1}^k m_{ij}}}{(\sum_{j=1}^k m_{ij})!}\prod_{j=1}^k \left( \frac{\Delta_{l_j}}{T} \right)^{m_{ij}}\,,
\end{align*}
where $T = \Lambda(S) = \sum_{i=1}^\infty \Delta_i \delta_{X_i}(S) = \sum_{i=1}^\infty \Delta_i$.  Readjusting the outermost product in the above equation, we get,
\begin{equation*}
\mathbb{P}(N_1,\dotsc, N_n| \Lambda) = \frac{\exp(-nT)}{\prod_{i=1}^n(\sum_{j=1}^k m_{ij})!}\prod_{j=1}^k \Delta_{l_j}^{\sum_{i=1}^n m_{ij}}\,.
\end{equation*}
Since, we are not interested in the actual points $X_{l_i}$'s, but only the number of occurrences of the different points in the point processes, that is, $[m_{ij}]_{(n,k)}$,  we can sum over every $k$-tuple of distinct atoms in the random measure $\Lambda$. Hence,  
\begin{align*}
&P([m_{ij}]_{(n,k)}| \Lambda) \\
&= \frac{\exp(-nT)}{\prod_{i=1}^n(\sum_{j=1}^k m_{ij})!}\sum_{\Delta_{l_1}\neq \Delta_{l_2}\neq\dots\neq\Delta_{l_k}}\prod_{j=1}^k \Delta_{l_j}^{\sum_{i=1}^n m_{ij}}\,,
\end{align*}
where the sum is over all subsets of length $k$ of the set $\{\Delta_1, \Delta_2, \dots\}$.
Finally, in order to compute the result, we need to take expectation with respect to the distribution of $\Lambda$. Towards that end, we note that only the weights of $\Lambda$ appear in the above equation. From section~\ref{section:CRM}, we know that the weights of a CRM with Poisson intensity measure $\rho(\dz)\mu(\dx)$form a Poisson process with mean measure $\mu(S)\rho(\dz)$. Hence, it is enough to take the expectation with respect to the Poisson process.
\begin{align}\label{conditional}
&P([m_{ij}]_{(n,k)}) \\
&= \frac{1}{\prod_{i=1}^n(\sum_{j=1}^k m_{ij})!}\mathbb{E}\Bigg[\exp{(-nT)} \\
& \left. \sum_{\Delta_{l_1}\neq \Delta_{l_2}\neq\dots\neq\Delta_{l_k}}\prod_{j=1}^k \Delta_{l_j}^{\sum_{i=1}^n m_{ij}}\right]\,,
\end{align}

The expectation can further be simplified by applying Proposition $2.1$ of~\cite{james2005bayesian}.
\begin{proposition}[\cite{james2005bayesian}]\label{BayesianCalculus}
Let $\mathcal{N}$ be the space of all $\sigma-$finite counting measures on $\mathbb{R}^+$, equipped with an appropriate $\sigma$-field. Let $f:\mathbb{R}^+\to \mathbb{R}^+$ and $g:\mathcal{N}\to \mathbb{R}^+$ be measurable with respect to their $\sigma$-fields. Then, for a Poisson process $N$ with mean measure $\mathbb{E}[N(\dx)] = \rho(\dx)$,
\begin{equation*}
\mathbb{E}\left[g(N)e^{-\sum_{\Delta \in N} f(\Delta)}\right] = \mathbb{E}\left[e^{-\sum_{\Delta \in N} f(\Delta)}\right] \mathbb{E}[g(\bar{N})]\,
\end{equation*}
where $\bar{N}$ is a Poisson process with mean measure $\mathbb{E}[N(\dx)] = e^{-f(x)}\rho(\dx)$.
\end{proposition}
Applying the above proposition to~\eqref{conditional}, we get
\begin{align*}
&P([m_{ij}]_{(n,k)}) = \frac{\mathbb{E}\left[e^{-\sum_{i=1}^\infty n\Delta_i}\right]}{\prod_{i=1}^n(\sum_{j=1}^k m_{ij})!} \\
&\times \mathbb{E}\left[\sum_{\Delta_{l_1}\neq \Delta_{l_2}\neq\dots\neq\Delta_{l_k} \in \bar{N}}\prod_{j=1}^k \Delta_{l_j}^{\sum_{i=1}^n m_{ij}}\right]
\end{align*}
where $\bar{N}$ is a Poisson process with mean measure $\mathbb{E}[N(\dz)] = e^{-nz}\rho(\dz)\theta $. The first expectation can be evaluated using Campbell's proposition and is given by $\exp\left(-\theta \int_{\mathbb{R}^+} (1-e^{-nz})\rho(\dz)\right)$. In order to evaluate the second expectation, we construct a new point process from $N^*$ on ${\mathbb{R}^+}^k$ by concatenating every set of $k$ distinct points in $\bar{N}$. The expression in the second expectation can then be rewritten as 
\begin{equation*}
\sum_{(\Delta_{l_1},\dotsc, \Delta_{l_k})\in N^*}\prod_{j=1}^k \Delta_{l_j}^{\sum_{i=1}^n m_{ij}}
\end{equation*}
By Campbell's proposition for point processes,
\begin{equation*}
\mathbb{E}\left[\sum_{\Delta \in N} f(\Delta)\right] = \int_{z \in \mathbb{R}^+}f(z)\rho(\dz)\,,
\end{equation*}
where $\rho(\dz) = \mathbb{E}[N(\dz)]$. Moreover, since the point process $N^*$ is obtained by concatenating distinct points in $N$, $\mathbb{E}[N^*(\dz_1,\dotsc, \dz_k)] = \prod_{j=1}^k \mathbb{E}[\bar{N}(\dz_j)] = \prod_{j=1}^k \theta e^{-nz}\rho(\dz_j)$, whenever $z_j$'s are distinct. Hence,

\begin{align*}
 &\mathbb{E}\left[\sum_{\Delta_{l_1}\neq \Delta_{l_2}\neq\dots\neq\Delta_{l_k} \in \bar{N}}\prod_{j=1}^k \Delta_{l_j}^{\sum_{i=1}^n m_{ij}}\right] \\
 &= \prod_{j=1}^k \int_{z \in \mathbb{R}^+} \theta  e^{-nz}z^{\sum_{i=1}^n m_{ij}}\rho(\dz)\,.
\end{align*}
Hence, the final expression for the marginal distribution of the set of counts for each latent feature is given by

\begin{align*}
P([m_{ij}]_{(n,k)}) =& \frac{\exp\left(-\theta \int_{\mathbb{R}^+}  (1-e^{-nz})\rho(\dz)\right) }{\prod_{i=1}^n(\sum_{j=1}^k m_{ij})!} \\
&\times \prod_{j=1}^k \int_{z \in \mathbb{R}^+} \theta  e^{-nz}z^{\sum_{i=1}^n m_{ij}}\rho(\dz)\,
\end{align*}
The above expression can be simplified by letting $\psi(t) =  \theta\int_{\mathbb{R}^+} (1-e^{-tz})\rho(\dz)$. Hence, $\psi^{(l)}(t) = (-1)^{l-1}\int_{\mathbb{R}^+} \theta  e^{-tz}z^{l}\rho(\dz)$. Hence, the above expression can be rewritten as
\begin{align*}
P([m_{ij}]_{(n,k)}) =& {(-1)^{\sum_{i=1}^n \sum_{j=1}^k m_{ij}-k}}\frac{\theta^k e^{-\theta\psi(n)} }{\prod_{i=1}^n(\sum_{j=1}^k m_{ij})!}\\
& \times \prod_{j=1}^k \psi^{(\sum_{i=1}^n m_{ij})}(n)
\end{align*}
\end{proof}

\textbf{Proof of Corollary~\ref{CRMPoissonCon}}
\begin{proof}
From proposition~\ref{CRMdistinct}, the distinct points in the point processes $N_i, 1\le i \le n$, form a Poisson process with mean measure $\frac{\mu(\dx)}{\mu(S)}\psi(n)$. Hence, the total number of distinct points $k$ is distributed as Poisson$(\psi(n))$. Hence, conditioning equation~\eqref{ESPF} with respect to $k$, we get the desired result.
\end{proof}

\textbf{Proof of Proposition~\ref{CRMnonAtomic}}
\begin{proof}
Let $N = \sum_{i=1}^n N_i$. From the arguments of proposition~\ref{CRMdistinct}, $N$ is a CRM, and hence, can be written as $N(\dx) = \int_{\mathbb{R}^+}z\bar{N}(\dz, \dx)$ for some Poisson process $\bar{N}$. Let $\Pi$ be the random collection of points corrsponding to $\bar{N}$. Now define a map $f:\mathbb{R}^+\times S \to S$ as the projection map on $S$, that is, $f(x,y) = y$ and $M = f(\Pi) = \{\{f(x,y): (x,y)\in \Pi\}\}$, where the double brackets indicate that $M$ is a multiset. The rest of the arguments remain same as in proposition~\ref{CRMdistinct} and proposition~\ref{CRMPoisson}.
\end{proof}

\textbf{Proof of Lemma~\ref{tables}}
\begin{proof}
Using Proposition~\ref{CRMnonAtomic} to marginalize $\Lambda_i$ from~\ref{CRMPP1}, we get that $[m_{ij}]_{1\le j \le r_{i\LargerCdot}}$ is distributed as CRM-Poisson$(\Phi(S),\rho, 1)$, that is,
\begin{align}
&P([m_{ij}]_{1 \le j \le r_{i\LargerCdot}}| \Phi(S)) \notag\\
&= \frac{\exp\left(-\Phi(S) \int_{\mathbb{R}^+} (1-e^{-z})\bar{\rho}(\dz)\right) }{(\sum_{j=1}^{r_{i\SmallerCdot}} m_{ij})!} \notag\\
&\times\prod_{j=1}^{r_{i\SmallerCdot}} \int_{z \in \mathbb{R}^+} \Phi(S)  e^{-z}z^{m_{ij}}\bar{\rho}(\dz)\,
\end{align}
Let $m_{i\LargerCdot} = \sum_{j=1}^{r_{i\LargerCdot}} m_{ij}$. Taking expectation with respect to $\Phi(S)$, we get the marginal distribution of  $[m_{ij}]_{1\le j \le r_{i\LargerCdot}}$, where $r_{i\LargerCdot}$ is also random. 
\begin{align}
&P([m_{ij}]_{1 \le j \le r_{i\LargerCdot}}) \notag\\
&= \mathbb{E}\left[\exp\left(-\Phi(S) \int_{\mathbb{R}^+} (1-e^{-z})\bar{\rho}(\dz)\right) \Phi(S)^{r_{i\LargerCdot}} \right] \notag\\
& \hspace{1cm}\times \frac{\prod_{j=1}^{r_{i\SmallerCdot}} \int_{z \in \mathbb{R}^+}  e^{-z}z^{m_{ij}}\bar{\rho}(\dz)}{ m_{i\LargerCdot}!}\, \label{eqn1}
\end{align}
It is given that $$h(u) = \mathbb{E}[e^{-u\Phi(S)}] $$ $$\bar{\psi}(u) = \int_{\mathbb{R}^+}(1-e^{-uz})\bar{\rho}(\dz)\, ,$$ 
Hence $$\frac{d^{r_{i\SmallerCdot}}}{du^{r_{i\SmallerCdot}}} h(u) = (-1)^{r_{i\LargerCdot}} \mathbb{E}\left[\Phi(S)^{r_{i\LargerCdot}} e^{-u\Phi(S)}\right]$$
$$\bar{\psi}^{(m_{ij})}(u)  = (-1)^{m_{ij}-1} \int_{\mathbb{R}^+}e^{-uz}z^{m_{ij}} \bar{\rho}(\dz)$$
Using the above results with $u=\bar{\psi}(1)$, equation~\eqref{eqn1} can be rewritten as 
\begin{align}
&P([m_{ij}]_{1 \le j \le r_{i\LargerCdot}}) \notag \\ 
&= (-1)^{m_{i\LargerCdot}} h^{(r_i\SmallerCdot)}(\bar{\psi}(1)) \frac{\prod_{j=1}^{r_{i\SmallerCdot}} \bar{\psi}^{(m_{ij})}(1)}{ m_{i\SmallerCdot}!}
\end{align}
\end{proof}
\end{document}